\magnification=\magstep1
\hsize=16 true cm
\input amssym.def
\input amssym.tex
\phantom{}
\vskip 3 true cm
\centerline{\bf SPHERICAL STEIN MANIFOLDS AND THE WEYL INVOLUTION}
\bigskip
\bigskip
\centerline {\bf Dmitri Akhiezer}
\bigskip
\bigskip
 \footnote{}{\smallskip \noindent The author was supported 
by the Russian Foundation for Basic Research, Grant
07-01-00230, and by SFB/TR12 "Symmetrien und Universalit\"at 
in mesoskopischen 
Systemen" of the Deutsche Forschungsgemeinschaft.}

\bigskip
\bigskip
\centerline{\bf 1. Introduction}
\bigskip
\bigskip
\noindent
Let $G$ be a connected reductive algebraic group over $\Bbb C$.
A normal algebraic $G$-variety $X$ is called
spherical if a Borel subgroup $B\subset G$ has an open orbit in $X$.
For $X$ affine this is the case if and only if 
the function algebra ${\Bbb C}[X]$ is a multiplicity free
$G$-module, i.e.,
if and only if 
every irreducible $G$-module appears in ${\Bbb C}[X]$ at most once;
see [Se], [VK].
 
The notion of a spherical variety
has been carried over to the category of complex spaces.
Namely, let $X$ be a normal complex space
and $K$ a connected 
compact Lie group acting on $X$ by holomorphic transformations.
Every element of the complexified Lie algebra $\goth g= {\goth k}\otimes 
\Bbb C$ gives rise to a holomorphic vector field on $X$.
The complex $K$-space 
$X$ is called spherical if there exists a point $x\in X$ such that
the holomorphic
tangent space $T_x(X)$ is generated by the vector fields from a Borel
subalgebra $\goth b \subset \goth g$.
Similarly to the algebraic setting, a normal Stein $K$-space
is spherical if and only if the Fr\'echet $K$-module
${\cal O}(X)$ is multiplicity free, i.e., every irreducible $K$-module
either occurs in ${\cal O}(X)$ with multiplicity 1 or does not occur at all.
Moreover, in that case $X$ is a $K$-invariant domain in a spherical affine
$K^{\Bbb C}$-variety; see [AH].

In order to state our main results,
we recall some terminology. Given a group $\Gamma $ acting
on two sets $X, \, Y$  and an automorphism $\alpha$ of $\Gamma $,
we say that a map $\mu : X \to Y$ is $\alpha $-equivariant
if
$\mu (\gamma x) = \alpha (\gamma )\mu (x) $ for all $x \in X,\,
\gamma \in \Gamma $.
A self-map $\mu : X \to X$ is called an involution if $\mu ^2 = {\rm Id}$.
By an involution of a group we will always mean
an involution which is a group automorphism.

In this paper we consider only non-singular Stein $K$-spaces.
Our goal is a characterization of spherical $K$-manifolds
in terms of $K$-equivariant antiholomorphic involutions. 

\bigskip
\noindent
{\bf Theorem 1.1}\ \it Let $X$ be a 
Stein manifold on which a connected compact Lie group
$K$ acts by holomorphic transformations. 
Then $X$ is spherical if and only if
there exists an antiholomorphic involution $\mu : X \to X$ such that 
$$\mu (x) \in Kx \eqno {(1)}$$ 
for all $x \in X$.

\bigskip
\noindent
\rm
According to a theorem
of J.Faraut and E.G.F.Thomas,  
the existence of $\mu $ with the above properties is sufficient for
$X$ to be spherical; see [FT], Theorem 3. 
In fact, the setting in [FT] is more general,
namely,  
the manifold $X$ need not be Stein and the group $K$ is not 
necessarily compact.
The conclusion is that, in the presence of
$\mu $ satisfying (1), 
any irreducible unitary 
$K$-module occurs in ${\cal O}(X)$ with multiplicity $\le 1$.
A simplified proof in our context is found in [AP]. 
The
converse statement is known 
in one special case, namely, if $X$ is an affine variety
homogeneous under the complexified group $G = K^{\Bbb C}$; see [AP].

We recall the definition of a Weyl involution
in Sect.2.
Let $\theta $ be a Weyl involution of $K$. 
The crucial role in our considerations is played by 
$\theta $-equivariant antiholomorphic
involutions of Stein manifolds acted on by $K$. 

\bigskip
\noindent
{\bf Theorem 1.2}\ \it Let 
$K$ be a connected compact Lie group,
$\theta : K \to K$ a Weyl involution, and $X$
a Stein manifold on which $K$ acts by holomorphic transformations. 
Assume that $X$ is spherical with respect to $K$.
Then there exists
a $\theta $-equivariant antiholomorphic involution $\mu : X \to X$.

\bigskip
\noindent
\rm 
In the paper of E.B.Vinberg and the author,
the Weyl involution appears 
as a tool bringing the notion of
weak symmetry into the theory of spherical varieties.
If $X$ is homogeneous under $G = K^{\Bbb C}$ then Theorem 1.2 
is a consequence of  
Theorem 5.11 in [AV]. At this stage,
the proof relies upon
the classification of reductive spherical subgroups in simple groups
given by Kr\"amer's list.
In the general case, $X$ embeds
as a $K$-invariant domain into 
the space of a vector bundle $G\times _H V$, where $H$ is a reductive
spherical subgroup in $G$ and $V$ is a spherical $H$-module
(see Sect.5). Therefore we need
an elementary version of Theorem 1.2 for
a $K$-module
$X = V$, which we will then apply
to the fiber. As long as $K$ is connected,
this result does not require sphericity of the $K^{\Bbb C}$-module
$V$. Moreover, 
Theorem 3.1 shows that there exists a $\theta $-equivariant 
involution $\mu : V \to V$, which is in fact antilinear,
and that such a $\mu $ is essentially unique if $V$
is irreducible. 
However, 
in order to accomplish the proof of Theorem 1.2, 
one
must be able to apply Theorem 3.1 to possibly disconnected groups 
acting on the fiber $V$.  
Our argument is based on
Theorem 4.3. Here is its simplified version, 
which might be of independent interest.

\bigskip
\noindent
{\it Let $G/H$ be a spherical homogeneous space,
$\sigma : G \to G$ an involution defining a split real form of $G$,
and $U$ a Lie subgroup of $G$ containing $H$,
such that $H \vartriangleleft U$ and $U/H$
is compact. If $\sigma (H) = H$ then
$\sigma (U) = U$ and the action of $\sigma $
on the (abelian) group $U/H$ is the inversion.}
\bigskip
\noindent
It is interesting to compare this statement with
Proposition 5.2 in [AV]
(see also Corollary 4.4).
The proof of main results is given in Sect.5. 

Throughout the paper,  
we denote by $N_G(H)$ the normalizer of a subgroup $H$
in a group $G$, by $G^\circ $ the connected component 
of the neutral element $e$ in a topological
group $G$, and by ${\rm X} (G)$
the character group of an algebraic group $G$.  
 
\bigskip
\bigskip
\centerline {\bf 2. Preliminaries}
\bigskip
\bigskip
\noindent
Let $G$ be a connected reductive algebraic group over $\Bbb C $.
Recall that an involutive regular automorphism $\theta : G \to G$
is called a Weyl involution if there exists a maximal
torus $T \subset G$ such that $\theta (t) = t^{-1}$ for $t\in T$. 
It is known that such an involution exists and that any two
Weyl involutions of $G$ are conjugate by an inner automorphism;
see e.g. [He], Ch. IX, $\S$5. 
From the point of view of the representation theory, the main
property of a Weyl involution is the following one. Let $V$ be a 
rational $G$-module,
$\varrho : G \to {\rm GL}(V)$ the corresponding representation, and
$V^\theta $ the $G$-module given
by the representation $g\mapsto \varrho (\theta (g))$ 
in the same vector space $V$. Then $V^\theta $ is isomorphic
to the dual $G$-module $V^*$.

Let $G^\theta  \subset G$ be the fixed point subgroup of $\theta $. 
This is a symmetric algebraic subgroup of $G$, which is always reductive. 
One can find a Cartan involution $\tau $ commuting with $\theta $ or, 
equivalently, a maximal compact subgroup $K\subset G$
stable under $\theta $. 
Any two maximal compact subgroups with this property are conjugate 
by an element of $G^\theta $.
Moreover, one can arrange that such a subgroup contains 
the maximal compact subgroup $T_c \subset T$; 
see [He], Ch. VI, A8 or Lemma 2.2 below. 
Thus, for any choice of $K$ satisfying 
$\theta (K) = K$ there is an element $a\in G^\theta $ such that
$K$ contains the compact torus $aT_ca^{-1}$. Since
$\sigma = \tau \theta $ preserves this torus and acts on it as the inversion, 
the fixed point subgroup of $\sigma $ is a split real form of $G$.

In what follows,
we will keep the above notation.
Namely, we will consider three involutive automorphisms of $G$. These are
a Weyl involution $\theta $, a Cartan involution $\tau $
commuting with $\theta $ and the product $\sigma = \theta \tau = \tau \theta $. The fixed point subgroups
$K =  G^\tau $ and $G^\sigma $
are a compact real form and, respectively, a split real form of $G$.
We remark that the notion of a 
Weyl involution can also be defined for connected compact Lie groups.
The definition is similar and its meaning 
for the representation theory is similar as well.
Note that the restriction $\theta \vert _K$ is a Weyl involution of $K$.  

For the convenience of the reader we recall the proof of two simple facts
about Lie groups
with finitely many connected components which we will need in the sequel.  
\bigskip
\noindent
{\bf Lemma 2.1}\ \it 
Let 
$\pi : G_1 \to G_2$ 
be an epimorphism of Lie groups. If ${\rm Ker} \, \pi $
has finitely many connected components and
$K_2 \subset G_2$ is a compact subgroup then there exists
a compact subgroup $K_1 \subset G_1$ such that $K_2 = \pi (K_1) $.
\bigskip
\noindent
Proof\ \rm 
Without loss of generality assume that $G_2 = K_2$.
Then $ G_1$ has finitely many connected components.
We have $K_2 ^\circ  = \pi (G_1^\circ ) $. 
By a theorem of D.Montgomery there exists a compact subgroup $L
\subset G_1^\circ $
such that $\pi (L) = K_2 ^\circ $; see [M], Cor. 2.
Let $K_1$ be a maximal compact subgroup of $G_1$ such that $L\subset K_1$.
Then $G_1 = K_1\cdot G_1^\circ $ by a theorem of G.D.Mostow; see [Mo].
Hence $K_2 = \pi (G_1) = \pi (K_1) \cdot K_2^\circ = \pi (K_1)$. 
 \hfill $\square $
\bigskip
\noindent
{\bf Lemma 2.2}\ \it Let $G_1\subset G_2$ be two Lie groups with 
finitely many connected components,
$G_1 $ a closed subgroup of $G_2$, 
and $\Delta $ a finite group of automorphisms of $G_2$, which leaves
invariant both $G_1$ and a maximal compact subgroup $K_1 \subset G_1$.
Then there exists a maximal compact subgroup $K_2 \subset G_2$, 
containing $K_1$
and invariant under each automorphism
from $\Delta $.

\bigskip
\noindent
Proof \ \rm
Consider the semidirect product $\tilde {G_2} = G_2 \rtimes \Delta $,
where $G_2$ is normal 
in $\tilde {G_2}$, with the 
given action of $\Delta $ on $G_2$ as an automorphism group. 
We identify $\Delta $ with the corresponding subgroup in $\tilde{G_2}$.
The semidirect products $\tilde {G_1} := G_1 \rtimes \Delta $ and 
$\tilde {K_1} := K_1\rtimes \Delta $
are subgroups of $\tilde {G_2}$.

Maximal compacts subgroups of a Lie group with finitely many connected
components are characterized 
by the two properties (see [Mo]): 
(i) such a subgroup has non-empty intersection
with every connected 
component of the ambient group; 
(ii) its intersection with the connected component
of the identity is a maximal compact subgroup in that component.
It follows from this criterion 
that $\tilde{K_1}$ is a maximal compact subgroup of $\tilde {G_1}$.
Let $\tilde{K_2}$ 
be a maximal compact subgroup of $\tilde{G_2}$ containing $\tilde{K_1}$.  
Since $\tilde {K_2}$ 
has non-empty intersection with all connected components of $\tilde {G_2}$
including those which are contained in $G_2$ 
(and are therefore connected components of $G_2$)
and, also, $(\tilde{G_2})^\circ = G_2^\circ $, 
we see from the above criterion that $K_2
:= \tilde {K_2}\cap G_2$ is a maximal compact subgroup of $G_2$. 
Because $K_1 \subset \tilde{K_1}
\subset \tilde{K_2}$, we have $K_1\subset K_2$. 
Furthermore, $\Delta \subset \tilde{K_1}\subset \tilde{K_2}$.
Thus, if $k\in K_2,\, \delta \in \Delta $ then 
$\delta k \delta ^{-1} \in \tilde { K_2} $.
Since $\delta G_2\delta ^{-1} = G_2$,  
we also have $\delta k \delta ^{-1} \in G_2$,
showing that $K_2$ is $\delta $-invariant.

\hfill $\square $  
\bigskip
\bigskip
\centerline{\bf 3. Antilinear $\theta $-equivariant involutions of 
$K$-modules}
\bigskip
\bigskip
\noindent
For a complex vector space $V$ let
$\overline {V}$  
denote
the conjugate vector space. By definition, 
$\overline {V} $ coincides with $ V$ as a real vector space, whereas 
the multiplication by complex
numbers on $\overline{V}$ is given by $(c,v) \mapsto \bar c v\ \  
(c\in \Bbb C,\ v\in V)$. Let $\iota : V \to \overline{V}$ 
denote the identity
map of real vector spaces, 
which is antilinear with respect to the complex structures on 
$V$ and $\overline{V}$. 
If $\varrho : K \to {\rm GL}(V)$ is a representation of a group $K$
then there is 
a natural representation of $K$ in $\overline{V}$ 
which we denote by $\bar \varrho $.
We have $\iota (\varrho (k) v) =\bar \varrho (k)  \iota (v)$ 
for all $ k\in K,\ v\in V$.
\bigskip
\bigskip
\noindent
{\bf Theorem 3.1}\ \it Let 
$K$ be a connected compact Lie group, 
$\theta : K \to K$ a Weyl involution, and $\varrho : K \to {\rm GL}(V)$
a complex 
representation of $K$. 
Then there exists an antilinear mapping 
$\mu : V \to V$, such that $\mu ^2 = {\rm Id}$
and
$$\mu (\varrho (k) v) = \varrho (\theta (k)) \mu (v) $$
for all $k \in K,\ v\in V$. If $\varrho $ is irreducible and 
$\mu ^\prime : V \to V$ is another mapping with the above properties then
$\mu^\prime = c\mu $ for some $c\in {\Bbb C}, \vert c \vert = 1$. 
\bigskip
\noindent Proof \ \rm A $K$-invariant 
positive Hermitian form on $V$ gives rise to an isomorphism
of $K$-modules $\overline{V} \simeq V^*$. 
Let $V^\theta $ be the $K$-module defined 
by the representation $k \mapsto \varrho (\theta (k))$ in $V$.
Then, by the main property of a Weyl involution, 
the $K$-modules $V^\theta $ and $\overline{V}$ are also isomorphic.
Thus we have an isomorphism of complex vector spaces 
$A : V \to \overline{V}$ such that
$$ A\, (\varrho (k) v )= \bar {\varrho}(\theta (k)) Av$$
for all $k\in K, v\in V$. Define the map $\mu : V \to V$ by 
$\mu = \iota ^{-1}\cdot A$. Then $\mu $
is antilinear and 
$\mu (\varrho (k) v)= \varrho (\theta (k))  \mu (v)$.

Note that it is enough to prove the theorem for irreducible 
representations. Indeed, if $V$ is the direct sum of irreducible $K$-modules
$V_j$ and $\mu _j : V_j \to V_j$ has the required property for each $j$ then 
$\mu $ can be defined as the direct sum of
$\mu_j$. We assume now that $V$ is an irreducible $K$-module.
Let $T\subset K$ be a maximal torus on which $\theta $ acts 
as the inversion and let $v\in V$ be a weight vector, such that
the corresponding weight $\lambda : T \to {\Bbb C}^* $ has multiplicity one. 
For the map $\mu $ constructed above and for any $t\in T$
we have
$$\varrho(t)^{-1} \mu(v) = \varrho(\theta(t)) \mu(v) 
= \mu(\varrho(t) v) 
=\mu (\lambda(t)v)=\overline {\lambda (t)}\mu (v)$$
or, equivalently,
$$ \varrho (t) \mu (v) = \lambda (t)\mu (v).$$
Therefore $\mu (v) = av$, where $a\in {\Bbb C}^*$.
Applying $\mu $ again we get $\mu ^2(v) = \vert a \vert ^ 2 v$.
Since $V$ is an irreducible $K$-module and $\mu ^2$ 
is a (complex linear) intertwining operator,
it follows by Schur's Lemma that $\mu ^2 = \vert a \vert ^2 \cdot {\rm Id}$.
Replacing $\mu $ by $\vert a \vert ^{-1}\mu $ we obtain the mapping with 
all required properties.
Finally, assume $\mu^\prime $ and $\mu $ are two such mappings.
Then $\mu (v) = av$ and $\mu^\prime (v) = a^\prime v$ as above. 
Because both $\mu ^2 $ and ${\mu ^ \prime }^2 $ are the identity maps,
it follows that $\vert a \vert = \vert a^\prime \vert = 1$.
Since $\mu \mu ^\prime $ is a (complex linear) intertwining 
operator and $\mu \mu ^\prime (v) = a(a^\prime )^{-1}v$,
we have $\mu \mu ^\prime = a(a^\prime)^{-1}\cdot {\rm Id}$. 
Multiplying this equality by $\mu $ on the left,
we get $\mu^\prime = c\mu $, where $c = a^\prime  a^{-1},\   
\vert c\vert = 1$.
\hfill {$\square $}  
\bigskip
\bigskip
\centerline{\bf 4. Involution $\sigma $ and spherical subgroups }
\bigskip
\bigskip
\noindent
Let $H\subset G$ be an algebraic subgroup, 
$\phi \in {\rm X}(H)$, and
${\Bbb L}^\phi $ the  
homogeneous line bundle over $G/H$ 
associated with the character $\phi $.
Recall that ${\Bbb L}^\phi = G\times _H {\Bbb C}$, where the action of $H$
on $G\times {\Bbb C}$ is given by
$$(g, z) \ {\buildrel h\over \mapsto} \ 
(gh^{-1}, \phi (h)z)\qquad (g\in G,  h\in H, 
z \in {\Bbb C}).$$
The
space of regular sections 
$\Gamma (G/H, {\Bbb L}^\phi )$ is a $G$-module. Namely,
the elements
of $\Gamma (G/H, {\Bbb L}^\phi )$ 
are identified with the regular functions on $G$
satisfying the equation
$f(gh) = \phi (h)^{-1}f(g)$, where $g\in G,\, h\in H$.
Under this identification, $G$
acts in $\Gamma (G/H, {\Bbb L}^\phi )$ by left translations. 
We write
$(g\cdot f)(x) = f(g^{-1}x)$ for $g,x\in G$ 
and $f\in \Gamma (G/H,{\Bbb L}^\phi )$. 

Let
$N = N_G(H)$ and assume that 
$$\phi (uhu^{-1}) = \phi (h) \qquad {\rm for}\ h\in H,\, u\in N. \eqno{(2)}$$
Since ${\rm X}(H)$ is a discrete group, 
this is automatic for $u \in HN^\circ $. Thus the orbit of $\phi $
under $N$ acting by conjugation is finite. 
This will be used in the proof of Theorem 4.3.
It follows from (2) that ${\rm Ker} \, \phi $ is normal in $N$. Also, 
$\Gamma (G/H, {\Bbb L}^\phi )$ has a structure of an $N$-module with the
action given by 
$$f\mapsto f^u,\ f^u(g) = f(gu),\, u\in N.$$ 
Indeed,
if $f$ is a section of ${\Bbb L}^\phi $ then
$$f^u(gh) = f(ghu) = f(gu\cdot u^{-1}hu) 
= \phi (u^{-1}hu)^{-1}f(gu) = \phi (h)^{-1} f^u(g),$$
where the last equality follows from (2).
This shows that $f^u $ is also a section of ${\Bbb L}^\phi $. 
We refer to this action as 
the right action of $N$ in $\Gamma (G/H, {\Bbb L}^\phi )$.

Recall that an algebraic subgroup $H \subset G$ is called spherical
if the homogeneous space $G/H$ is a spherical variety.
In what follows,
we are mainly interested in spherical subgroups $H \subset G$.
Though this fact is not used in the proofs below, it is good to remember
that for $H$ spherical the group $N/H$ is diagonalizable and,
in particular, abelian; see [BP].

\bigskip
\noindent
{\bf Theorem 4.1}\ \it Let $H \subset G$ be a spherical subgroup,
$a\in G$, and
$\sigma (H) = aHa^{-1} $. Let
$\phi \in {\rm X}(H) $ be a character satisfying (2). Assume in
addition that 
$$\phi (a^{-1}\sigma (h) a) = \overline {\phi (h)} \eqno {(3)}$$
for all $h\in H$.
Then any 
irreducible $G$-submodule $V \subset \Gamma (G/H, {\Bbb L}^\phi )$
is invariant under the right action of $N$. Moreover, $N$ acts on $V$ by 
$$ f^u = \phi _{_V}(u)^{-1}\cdot f  \ \ \qquad (f\in V, \, u\in N),$$
where the character
$\phi _{_V} \in {\rm X}(N)$ satisfies
$$\overline {\phi _{_V}(u)} = \phi _{_V}(\sigma (aua^{-1})) \eqno {(4)}$$
(this makes sense because $\sigma (N) = aNa^{-1}$).

\bigskip
\noindent
Proof \ \rm 
Fix a section $f\in \Gamma (G/H, {\Bbb L}^\phi )$ and put
$$F (g) = \overline {f(\sigma (g)a) } .$$
Clearly, $F$ is a regular function on $G$, and
$$F(gh) = \overline { f(\sigma (g)\sigma (h)a)} 
= \overline { f(\sigma (g)a \cdot a^{-1}\sigma(h)a)} =
$$
$$=
\overline {\phi (a^{-1}\sigma(h)a) ^{-1}} 
\cdot \overline {f(\sigma (g)a)} 
= \phi (h) ^{-1} F(g)$$
in view of (3).
Therefore $F\in \Gamma (G/H, {\Bbb L}^\phi)$, and so we obtain
the antilinear map
$$\eta : 
\Gamma (G/H, {\Bbb L}^\phi) 
\to  
\Gamma (G/H, {\Bbb L}^\phi), \ \qquad 
\eta (f)(g) = \overline { f(\sigma (g)a)}.$$
For $k\in K$ we have   
$$\eta (k\cdot f)(g) = \overline{ (k\cdot f) (\sigma(g)a)}
= \overline { f(k^{-1}\sigma(g)a)}
= \overline{f(\sigma (\theta (k)^{-1}g)a) } =
$$
$$
= \eta (f) (\theta (k)^{-1}g)
= [\theta (k)\cdot \eta (f)] (g)$$
because $\theta $ and $\sigma $ coincide on $K$.
Therefore $\eta (V)$ is a $K$-submodule and hence a $G$-submodule in
$\Gamma (G/H, {\Bbb L}^\phi )$. Moreover, the mapping $\eta : V \to \eta (V)$
gives rise to a 
$K$-equivariant antilinear bijective map $V \to \eta (V)^{\theta }$,
where the $K$-module $\eta (V)^{\theta }$ 
is obtained from $\eta (V)$ by changing the corresponding
representation $k \mapsto 
\varrho (k)$ to $k\mapsto \varrho (\theta (k)) \ (k\in K)$.
Let $\overline {V}$ be the complex conjugate space of $V$ with
the natural representation of $K$. As we have seen, there is an
isomorphism of $K$-modules $\overline {V}
\simeq \eta (V)^\theta $. 
On the other hand,
from the existence of a $K$-invariant positive Hermitian form on $V$
we get an isomorphism of $K$-modules
$\overline{V} \simeq V^*$. Furthermore, since $\theta $ is a Weyl involution, 
we also have the isomorphism
$\eta (V)^\theta \simeq \eta (V)^*$.  
Combining these isomorphisms,
we conclude that $V^* \simeq \eta (V)^*$ and $V\simeq \eta (V)$.
But $H$ is a spherical subgroup in $G$, and so each irreducible $G$-submodule
in $\Gamma (G/H, {\Bbb L}^\phi )$ has multiplicity one. Thus $\eta (V) = V$.

The action of $G$ 
commutes with the right action of $N$ on $\Gamma (G/H, {\Bbb L}^\phi ) $. 
Since 
$\Gamma (G/H, {\Bbb L}^\phi )$ 
is multiplicity free, $V$ is invariant under the
action of $N$. By Schur's Lemma $N$ acts on $V$ via multiplication by a 
character. We denote the latter one by 
$ \phi _{_V} (u) ^{-1} , \ u\in N$.
In other words, for $f\in V $ we have
$$f(gu) = \phi _{_V}(u)^{-1}f(g).$$
Replacing $g$ by $\sigma (g)a$, we get
$$f(\sigma(g)au) = \phi _{_V} (u)^{-1}f(\sigma (g)a) 
= \phi _{_V}(u)^{-1} \cdot \overline {\eta (f)(g)}.$$ 
On the other hand,
let $u^\prime = \sigma (aua^{-1})$. Then 
$$f(\sigma (g)au) = f(\sigma (g)aua^{-1}a) = f(\sigma (gu^\prime)a) 
= \overline { \eta (f)(gu^\prime )} =
\overline {\phi _{_V} (u^\prime ) ^{-1}} \cdot \overline {\eta (f)(g)},$$
where in the last equality we used that $\eta (f) \in V$. Together with
the previous computation, this proves (4).
\hfill $\square $
\bigskip
\noindent

\bigskip
\noindent
{\bf Corollary 4.2}\ \it Let $H\subset G$ 
be a spherical subgroup, such that $\sigma (H) = aHa^{-1}$
for some $a \in G$.
Assume that  
$\phi \in {\rm X}(H)$ 
satisfies (2) and (3) and let $U\subset N$ be
a Lie subgroup with compact image in $N/{\rm Ker}\, \phi $. 
Then for any $u\in U$ the right action of $\sigma (aua^{-1})$
on $\Gamma (G/H, {\Bbb L}^\phi )$ coincides with that of $u^{-1}$.
\bigskip
\noindent
Proof\ \rm 
Applying Lemma 2.1 to the epimorphism $N \to N/{\rm Ker}\, \phi $,
we can write any $u\in U$ as $u = ch$, where $h\in H, \, \phi (h) = 1$, and $c$
is contained in a compact subgroup of $N$. 
Let 
$V \subset \Bbb C[G/H]$ be an irreducible $G$-submodule. 
Since $\phi _{_V} \vert _H = \phi $ by the definition of $\phi _{_V}$,
we have
$\phi _{_V}(u) = \phi _{_V}(c)$ and, therefore,
$\vert \phi _{_V}(u) \vert = 1$. By Theorem 4.1, it follows from (4)
that 
$$\phi _{_V} (u^{-1}) = \phi _{_V}(\sigma (aua^{-1})).$$ 
This is true for any $G$-submodule $V \subset {\Bbb C}[G/H]$, showing that 
$u^{-1}$ and $\sigma (aua^{-1})$ induce the same linear
transform of $\Gamma (G/H, {\Bbb L}^\phi )$. 
\hfill $\square $ 
\bigskip
\noindent
{\bf Theorem 4.3}\ \it Let $H\subset G$ 
be a spherical subgroup, $a\in G$, and
$\sigma (H) = aHa^{-1}$.
Let $U\subset N$ be a Lie subgroup 
with compact image in $N/H$. Then
$$\sigma(aua^{-1}) \equiv u^{-1} \qquad ({\rm mod}\ H) $$
for all $u\in U$.
In particular, if $H \subset U$ then
$U$ is stable under the mapping $u \mapsto \sigma (aua^{-1})$
and the induced automorphism of $U/H$ is the inversion.
\bigskip
\noindent
Proof\ \rm
By Chevalley's theorem there exist a rational representation
$\rho : G \to {\rm GL}(V)$ and a vector $v \in V$, such that
$$H = \{g \in G \ \vert \ \rho (g) v\in {\Bbb C}\cdot v\}.$$
Let $\phi \in {\rm X}(H)$ denote the corresponding character,
i.e., $\rho (h)v = \phi (h)\cdot v$ for $h\in H$.
For $u \in N$ define $\phi ^u\in {\rm X}(H)$ by
$$\phi ^u(h) = \phi (uhu^{-1})\quad (h\in H).$$ 
Since 
the family $\{\phi ^u\}_{u\in N}$
is finite, we can find $u_1,\ u_2,\ \ldots ,\ u_k \in N$
with the property that for every $u\in N$ there exists exactly one
$u_i,\ i\in [1,k]$, such that $\phi ^u = \phi ^{u_i}$. Write $\phi _i$
instead of $\phi ^{u_i}$, define the vectors $v_i \in V$
by $v_i = \rho (u_i)v$, and consider the tensor product
of representations $\tilde \rho : G \to {\rm GL}(\tilde V)$, where
$\tilde V = V\otimes \ldots \otimes V$ ($k$ times).  
For the decomposable vector $\tilde v = v_1\otimes \ldots \otimes v_k \in 
\tilde V$
we have
$$H = \{ g \in G \ \vert \ \tilde \rho (g)\tilde v
\in {\Bbb C}\cdot \tilde v\},$$
and the corresponding character of $H$ is the product
$$\tilde \phi = \prod _{i=1}^k \ \phi _i,$$
which is invariant under $N$, i.e.,
$\tilde \phi (uhu^{-1}) = \tilde \phi (h)$ for all $h\in H, \ u\in N$.
Therefore, replacing 
the representation
$\rho : G \to {\rm GL}(V)$ and 
the vector $v\in V$ by
$\tilde \rho : G \to \tilde {\rm GL}(\tilde V)$ and,
respectively, $\tilde v \in \tilde V$
we may assume that $\phi $ satisfies (2).

We now want to change the triple $V, v$ and $\rho $ 
without changing $H$, so that $\phi $ also satisfies (3).
In the notation of Sect. 3, 
let $\overline V$ denote the complex conjugate vector space. 
We have the antilinear
map $\iota : V \to \overline V$ 
and the antiholomorphic representation
$\bar \rho : G \to {\rm GL}(\bar V)$ such that
$\iota (\rho (g) ) = \bar \rho (g) \iota $ for $g\in G$.
Define a representation of $G$ in the tensor product $V \otimes \overline V$ 
as $g\mapsto \rho (g)$ 
on the first factor and $g\mapsto \bar \rho (a^{-1}\sigma (g)a)
$ on the second one.
Note that this representation is rational. Moreover, it is easily seen
that $H$ is the stabilizer of the line 
${\Bbb C} \cdot v\otimes \iota (v)$
and the corresponding character is given by
$$\hat \phi ( h) =  \phi (h)\cdot \overline {\phi (a^{-1}\sigma (h)a)},\quad 
h \in H.$$ 
Since $\sigma (a)a\in N$ and $\phi ^u = \phi $
for $u\in N$, we see that
$\hat \phi $ satisfies (3).

We just proved that there exists an equivariant
embedding $j: G/H \to {\Bbb P}(V)$ with $j^*{\cal O}(-1) = {\Bbb L}^\phi $.
Thus the sections of the dual bundle $({\Bbb L }^\phi )^* = 
{\Bbb L}^{\varphi }$, where $\varphi = \phi ^{-1}$, separate
points of $G/H$. It follows that the 
ineffectivity kernel of the right action of
$N$ 
in $\Gamma (G/H , {¸\Bbb L}^\varphi )$
is exactly $H$.
If $U \subset N$ is a Lie subgroup with compact image in $N/H$,
then by Lemma 2.1 there exists a possibly larger Lie subgroup of $N$
with compact
image in $N/{\rm Ker}\ \varphi $, which is mapped onto $UH/H$.
We denote this larger subgroup again by $U$. Let $u \in U$.
Then, by Corollary 4.2, 
$u^{-1}$ 
and $\sigma (aua^{-1})$ induce the same transform of
$\Gamma (G/H, {\Bbb L}^\varphi )$.
Therefore $\sigma (aua^{-1}) \equiv u^{-1}$ modulo $H$.
\hfill $\square $ 

\bigskip
\noindent
As a corollary, we get another proof of Proposition 5.2 in [AV].
\bigskip
\noindent 
{\bf Corollary 4.4}\ \rm (see [AV], Proposition 5.2)
\it \ Let $H \subset G$ be a reductive spherical subgroup invariant
under a Weyl involution $\theta $. 
Then $\theta $ acts as the inversion on the abelian group $N/H$. 
\bigskip
\noindent
Proof\ \rm Take $\Delta = \{e,\, \theta\}$ in Lemma 2.2
and apply the assertion of the lemma to each of the three inclusions
$\{e\} \subset H \subset N \subset G$, moving from the left to the right.
Then we get $\theta $-stable 
maximal compact subgroups in $H,\, N$ and $G$ together
with similar inclusions. In other words, there is 
a $\theta $-stable maximal compact subgroup $K \subset G$, such that
$K\cap H $ and $K\cap N$ are 
maximal compact subgroups in $H$ and, respectively, in $N$.
Let $\tau $ 
be the Cartan involution of $G$ with $G^\tau = K$. 
Clearly, $\theta \tau = \tau \theta $.
Since $\tau $ fixes pointwise 
the maximal compact subgroups $K\cap H$ and $K\cap N$
of the reductive groups $H$ and $N$, it follows that 
$\tau (H) = H$ and $\tau (N) = N$. Thus 
$\tau $ acts on $N/H$. 
Note that $N/H$ 
is a reductive algebraic group and 
$(K\cap N) / (K \cap H)$ is its maximal compact
subgroup. 
By Theorem 4.3 the product $\sigma = \tau \theta  $ 
and therefore also $\theta $
acts as the inversion on $(K\cap N )/ (K\cap H)$,
showing that the latter group is abelian. 
(Of course, it follows that $N/H$ is abelian, but this fact is known in a 
more general setting; see [BP]).
Since $\theta $ is an 
algebraic automorphism, $\theta $ 
acts as the inversion on the whole group $N/H$.
\hfill $\square $   
\bigskip
\bigskip
\centerline {\bf 5. Proof of main results}
\bigskip
\bigskip
\noindent
\it Proof of Theorem 1.2 \ \rm
First of all, we reduce the statement to the algebraic case.
Namely, by Theorem 2 in [AH], $X$ can be embedded as a $K$-invariant 
domain into
an affine spherical $K^{\Bbb C}$-variety
which we denote by $X^{\Bbb C}$. Moreover, if $G = K^{\Bbb C}$
and $i: X \hookrightarrow X^{\Bbb C}$ is the embedding then 
$G\cdot i(X) = X^{\Bbb C}$. Hence $X^{\Bbb C}$ is non-singular.
Assuming that the theorem is
proved for $X^{\Bbb C}$, we get a $\theta $-equivariant antiholomorphic
involution $\mu : X^{\Bbb C} \to X^{\Bbb C}$. Since the variety $X^{\Bbb C}$
is spherical, its Fr\'echet algebra ${\cal O}(X)$ is a multiplicity
free $K-$module. Thus, by Theorem 4.1 in [AP], we have $\mu (x) \in Kx$ for
all $x\in X^{\Bbb C}$.
It follows that $\mu (i(X)) = i (X)$, and so the required involution of $X$
is obtained from $\mu $ by restriction.  

From 
now on we assume that $X$ is a smooth affine spherical variety of $G$.
Since $\Bbb C[X]={\Bbb C}$, 
the well-known application of Luna's Slice Theorem displays $X$
as a vector bundle. Namely,
$X=G\times _{H}V$, where $H$ is a reductive subgroup
of $G$ and $V$ is an $H$-module; see
[Lu], Cor.2, p.98. Moreover, since 
$X$ is a spherical variety, it follows that $H $ is a spherical
subgroup in $G$ and $V$ is a spherical $H$-module; see [K-VS], Corollary 2.2.

It suffices to prove our theorem for some maximal compact subgroup
$K \subset G$ and for some Weyl 
involution $\theta $ of $G$ satisfying $\theta (K) = K$.
Let $L$ be a maximal compact subgroup of $H$ contained in $K$. 
According to [AV], Prop. 5.14, we can find
a Weyl involution $\theta $ of $G$ so that
$K$ and $L$ are $\theta $-invariant and $\theta $
induces a Weyl involution of $K$ and $L^{\circ }$.
Then, of course,
$\theta (H) = H$
and $\theta \vert _{H^\circ}$ is a Weyl involution of $H^{\circ }$.
Moreover, if $\tau : G \to G$ is the Cartan involution with $G^\tau = K$
then $\tau \theta = \theta \tau $. We
will use our standard notation $\sigma = \tau \theta$. 

By Theorem 3.1
there exists an antilinear involution $\nu : V \to V$, such that
$\nu (lv) = \theta (l) \nu (v)$ for $v\in V$ and $l \in L^{\circ }$.
Since
$L^\circ $ is a maximal compact subgroup of
the connected reductive group $H^\circ $, it follows that
$\nu (hv) = \sigma (h)\, \nu (v)$
for $h \in H^\circ , \ v\in V$. Indeed, 
if $v$ is fixed then this equality holds for all $h \in L^\circ $.
But 
$L^\circ $ is a 
maximal totally real 
submanifold 
in
$H^\circ $, and so the above equality holds also on $H^\circ $. 
We claim that in fact
$$\nu (hv) = 
\sigma (h) \, \nu (v) \qquad {\rm for\ all} \ \ v \in V,\, h\in H.
\eqno {(*)}$$
\bigskip
\noindent
\it Proof of $(*)$.
\rm Take a vector $v_0 \in V$ such that the orbit
$G v_0$ is open in $X$. Then the orbit 
$H v_0$ is open in $V$.
Let $U = G_{v_0} = H_{v_0}$
be the isotropy subgroup at $v_0$. Put $v_1 = \nu (v_0)$.
Then $\nu (hv_0)= \sigma (h)v_1$ for $h \in H^\circ $,
hence $\nu (H^\circ v_0) = H^\circ v_1 $. 
But the open orbit $H v_0
\subset V$ is connected. Thus $H^\circ $ acts on this orbit transitively,
i.e., $Hv_0 = H^\circ v_0$. Since the orbit $H^\circ v_1$ is also open,
we have $H^\circ v_0 = H^\circ v_1$ by
the uniqueness of an open orbit.
Therefore $v_1 = av_0$, where $a\in H^\circ $.

Since $U \subset H$, the intersection
$U\cap H^{\circ }$ has finite index in $U$. Thus 
$G/(U\cap H^\circ) $ is  a spherical variety
together with $G/U$.
The equality $v_1 = av_0$ shows that the isotropy subgroup
of $v_1$ in $H^\circ $ is $a(U\cap H^\circ )a^{-1}$.
On the other hand,
$\nu $ is $\sigma $-equivariant with respect to the connected group
$H^\circ $. From this 
it follows that  
$a(U\cap H^\circ )a^{-1}= \sigma (U \cap H^{\circ })$. 
We now apply Theorem 4.3 to the pair $U \cap H^\circ \subset G$
and to the subgroup $U$ normalizing $U \cap H^\circ $. The quotient
$U/(U\cap H^\circ )$ is finite, 
and so we obtain
$\sigma (U) = aUa^{-1}$.

Since $H^\circ $ is transitive on $H v_0$, we have $H= H^\circ \cdot U$.
For any $h \in H$ write $h = h^\prime u$ with $h^\prime
 \in H^\circ,\, u \in U$.
Then
$$\nu (hv_0) = \nu (h^\prime uv_0) = \nu (h^\prime v_0) 
= \sigma (h^\prime )v_1$$
and, on the other hand,
$$\sigma (h)v_1 = \sigma (h^\prime )\sigma (u) v_1 
= \sigma (h^\prime ) au^\prime a^{-1}v_1 = \sigma (h^\prime )v_1 ,$$
where $\sigma (u) = au^{\prime }a^{-1}$ and $u^\prime \in U$. 
These computations show that 
$\nu (hv_0) = \sigma(h)v_1 $ 
for any $h \in H$.
Replacing $h$ by $hh^\prime $ where $h^\prime \in H^\circ $, we get
$\nu (hh^\prime v_0) = \sigma (hh^\prime )v_1 
= \sigma (h)\sigma (h^\prime )v_1 =
\sigma (h)\nu (h^\prime v_0)$. Hence
$\nu (hv) = \sigma (h) \nu (v)$
where $v = h^\prime v_0$ is any vector from the open orbit in $V$. 
Since the open orbit is dense, this completes the proof
of $(*)$.
\bigskip
\noindent
Let $\mu $ be the self-map of $G\times V$
defined
by
$$\mu (g,v) = (\sigma (g), \nu (v)).$$
Obviously, $\mu $ is an antiholomorphic involution of $G\times V$.
For $h\in H$ let $t_h (g,v) = (gh^{-1}, hv)$. 
This defines an action of $H$ on $G\times V$ and, 
by definition, $X$ is the geometric quotient with respect to this action. 
It follows from $(*)$ that $\mu \cdot t_h = t_{\sigma (h)}\cdot \mu $.
Thus $\mu $ defines a self-map of $X$ with all required properties. 
\hfill $\square $
\bigskip
\noindent
{\it Proof of Theorem 1.1}\ 
As we already noted in the introduction, the 
existence of an antiholomorphic involution $\mu : X \to X$ 
with property (1) implies that $X$ is spherical; see [FT], Theorem 3.
Conversely, if $X$ is a spherical Stein manifold then,
by Theorem 1.2, there exists
a $\theta $-equivariant antiholomorphic involution $\mu : X \to X$.
On the other hand, a holomorphically separable complex $K$-manifold
with a $\theta $-equivariant antiholomorphic involution is
multiplicity free if and only if (1) is fulfilled; see [AP], Theorem 4.1.
\hfill $\square $

\bigskip
\bigskip
\bigskip
\centerline {\bf References}

\bigskip
\bigskip
\noindent
[AV] D.N.Akhiezer, E.B.Vinberg, 
\it Weakly symmetric spaces and spherical varieties, 
\rm Transformation Groups 4 (1999), 3 - 24 
\smallskip
\noindent
[AH] D.Akhiezer, P.Heinzner,  
\it Spherical Stein spaces, \rm Manuscripta Math. 114 (2004), 327 - 334
\smallskip
\noindent
[AP] D.Akhiezer, A.P\"uttmann, \it 
Antiholomorphic involutions of spherical complex 
spaces, 
\rm Proc.Amer.Math.Soc., 
posted on January 3, 2008, PII S 0002-9939(08)08988-0 
(to appear in print) 
\smallskip
\noindent
[BP] M.Brion, F.Pauer, 
\it Valuations des espaces homog\`enes sph\'eriques, 
\rm Comment. Math. Helv. 62 (1987), 265 - 285
\smallskip
\noindent
[FT] J.Faraut, E.G.F.Thomas, 
\it Invariant Hilbert spaces of holomorphic functions, 
\rm J. of Lie Theory, 9 (1999), 383 - 402
\smallskip
\noindent
[He] S.Helgason, 
\it Differential Geometry, Lie Groups, and Symmetric Spaces, 
\rm Acad. Press, 1978.
\smallskip
\noindent
[K-VS] F.Knop, B. Van Steirteghem, 
\it Classification of smooth affine spherical varieties, \rm 
Transformation Groups, 11 (2006), 495 - 516
\smallskip
\noindent
[L] D.Luna, \it Slices \'etales, 
\rm Mem.Soc.Math.France (N.S.) 33 (1973), 81 - 105
\smallskip
\noindent
[M] D.Montgomery, \it Simply connected homogeneous spaces,
\rm Proc.Amer.Math.Soc. 1 (1950), 467 - 469
\smallskip
\noindent 
[Mo] G.D.Mostow, \it Self-adjoint groups, 
\rm Annals of Math. 62 (1955),
44 - 55 
\smallskip
\noindent
[Se] F.J.Servedio, \it Prehomogeneous vector spaces and varieties,
\rm Trans.Amer.Math.Soc. 176 (1973), 421 - 444
\smallskip
\noindent
[VK] E.B.Vinberg, B.N.Kimel'feld, \it Homogeneous domains in flag manifolds
and spherical subgroups of semisimple Lie groups,
\rm Funktsional'nyi Analiz i ego Prilozheniya 12 (1978), 12 - 19
(Russian); English translation: Funct.Anal.Appl. 12 (1978), 168 - 174 

\bigskip
\bigskip
\bigskip
\noindent Institute for Information Transmission Problems,
\smallskip
\noindent B.Karetny per. 19, Moscow, 101447, Russia

\end